\documentclass[reqno,centertags]{amsart}
\usepackage{amsmath}
\usepackage{amssymb}
\usepackage{amsfonts}

\newtheorem{theorem}{Theorem}
\theoremstyle{plain}

\numberwithin{equation}{section}

\input{tcilatex}

\begin{document}
\title[]{Surfaces of coordinate finite $II$-type}
\author{Hassan Al-Zoubi}
\address{Department of Mathematics, Al-Zaytoonah University of Jordan, P.O. Box 130, Amman, Jordan 11733}
\email{dr.hassanz@zuj.edu.jo}
\author{Mutaz Al-Sabbagh}
\address{Department of Basic Engineering Sciences, Imam Abdulrahman bin Faisal University, Dammam 31441, Saudi Arabia}
\email{malsbbagh@iau.edu.sa}
\author{Tareq Hamadneh}
\address{Department of Mathematics, Al-Zaytoonah University of Jordan, P.O.
Box 130, Amman, Jordan 11733}
\email{t.hamadneh@zuj.edu.jo}
\date{}

\begin{abstract}
In this article, we study the class of surfaces of revolution in the 3-dimensional Euclidean space $E^{3}$ with nonvanishing Gauss curvature whose position vector $\boldsymbol{x}$ satisfies the condition $\Delta^{II}\boldsymbol{x}=A\boldsymbol{x}$, where $A$ is a square matrix of order 3 and $\Delta^{II}$ denotes the Laplace operator of the second fundamental form $II$ of the surface. We show that a surface of revolution satisfying the preceding relation is a catenoid or part of a sphere.
\end{abstract}

\maketitle

\section{Introduction}
Surfaces of finite $I$-type is one of the main topics that attracted the interest of many differential geometers from the moment that B. Y. Chen introduced the notion of surfaces of finite $I$-type with respect to the first fundamental form $I$ about four decades ago. Many results concerning this subject has been collected in \cite{C7}.

Let $\boldsymbol{x}:M^{2}\rightarrow E^{3}$ be a parametric representation of a surface in the $3$-dimensional Euclidean space $E^{3}$. Denote by $\Delta ^{I}$ the second Laplace operator according to the first fundamental form $I$ of $M^{2}$ and by $\boldsymbol{H}$ the mean curvature field of $M^{2}$. Then, it is well known that \cite{C3}
\begin{equation}
\Delta ^{I}\boldsymbol{x}=-2\boldsymbol{H}.  \notag
\end{equation}

Moreover in \cite{T1} T.Takahashi showed that the position vector $\boldsymbol{x}$ of $M^{2}$ for which $\Delta^{I}\boldsymbol{x} =\lambda\boldsymbol{x}$, 
is either the minimals with eigenvalue $\lambda=0$ or $M^{2}$ lies in an ordinary sphere $S^{2}$ with a fixed nonzero eigenvalue.

O. Garay in his article \cite{G3} has made a generalization of T. Takahashi's condition. Actually, so far, he studied surfaces in $\mathbb{E}^{3}$ satisfying $\Delta^{I}r_{i} = \mu_{i}r_{i}, i = 1,2,3$, where $\left(r_{1},r_{2},r_{3}\right)$ are the coordinate functions of $\mathbf{r}$. Another general problem also was studied in \cite{D3} for which surfaces in $\mathbb{E}^{3}$ satisfying $\Delta^{I}\mathbf{r}= K\mathbf{r} + L (\S)$, where $K \in M (3\times3) ;L \in M(3\times1)$. It was proved that minimal surfaces, spheres, and circular cylinders are the only surfaces in $\mathbb{E}^{3}$ satisfying $(\S)$. Surfaces satisfying $(\S)$ are said to be of coordinate finite type.

In the framework of the theory of surfaces of finite $I$-type in $E^{3}$, a general study of the Gauss map was made within this context in \cite{C9}. So, one can ask which surfaces in $E^{3}$ are of finite $I$-type Gauss map. On the other hand, it is also interested to study surfaces of finite $I$-type in $E^{3}$ whose Gauss map $\boldsymbol{n}$ satisfies a condition of the form $\Delta^{I}\boldsymbol{n}= A\boldsymbol{n}$, where $A \in\mathbb{Re}^{3\times3}$. Surfaces in $E^{3}$ whose Gauss map is of coordinate finite type corresponding to the first fundamental form was investigated by many researchers as one can see in \cite{A8, A10, B1, B2, B5, B6, D4}. Applications of similar results for polynomial functions over different domains are also given in \cite{H2, H3}.

In 2003 S. Stamatakis and H. Al-Zoubi in \cite{S1} followed the ideas of B. Y. Chen, they introduced the notion of surfaces of finite type regarding to the second or third fundamental forms, and since then much work has been done in this context.

Actually, ruled surfaces \cite{A4}, tubes \cite{A5}, quadrics \cite{A6} and a special case of surfaces of revolution \cite{Ar} are the classes of surfaces were studied in terms of finite type classification with respect to the third fundamental form. Mean while, tubes \cite{A3} and ruled surfaces \cite{A7} are the only classes of surfaces were investigated in terms of finite type classification with respect to the second fundamental form.

Another generalization can be made by studying surfaces in $E^{3}$ whose position vector $\boldsymbol{x}$ satisfies the following condition
\begin{equation} \label{3}
\Delta ^{J}\boldsymbol{x} =A\boldsymbol{x},  \ \ \ J = II, III,
\end{equation}%
where $A\in \mathbb{Re}^{3\times 3}$.

Regarding the third fundamental form in the 3-dimensional Euclidean space, it was proved that spheres and catenoids are the only surfaces of revolution satisfying condition (\ref{3}) \cite{S2}. Next, In \cite{A1} authors found that helicoids are the only ruled surfaces that satisfy (\ref{3}), meanwhile, spheres are the only quadric surfaces that satisfy (\ref{3}). Finally, in \cite{A2} it was shown that Scherk's surface is the only translation surface that satisfies (\ref{3}).

\section{Basic concepts}
Let $M^{2}$ be a smooth surface in $E^{3}$ parametrized by $\boldsymbol{x} = \boldsymbol{x}(u^{1}, u^{2})$ on a region $U: = (a, b) \times \mathbb{R}$ whose Gaussian curvature never vanishes. 
The standard unit normal vector field $\boldsymbol{n} $ on $M^{2}$ is defined by
\begin{equation}  \label{4}
\boldsymbol{n}=\frac{\boldsymbol{x_{u^{1}}}\times \boldsymbol{x_{u^{2}}}}{\|\boldsymbol{x_{u^{1}}}\times \boldsymbol{x_{u^{2}}}\|},
\end{equation}
where $\boldsymbol{x_{u^{1}}} := \frac{\partial\boldsymbol{x}(u^{1}, u^{2})}{\partial u^{1}}$ and $"\times"$ denotes the Euclidean vector product. We denote by
\begin{equation}  \label{5}
I = g_{ij}du^{i}du^{j},\ \ \    II = b_{ij}du^{i}du^{j},
\end{equation}
the first and second fundamental forms of $M^{2}$ respectively, where we put
\begin{equation*}
g_{11}= <\boldsymbol{x_{u^{1}}}, \boldsymbol{x_{u^{1}}}>,\ \ \ g_{12}= <\boldsymbol{x_{u^{1}}}, \boldsymbol{x_{u^{2}}}>,\ \ \ g_{22}= <\boldsymbol{x_{u^{2}}}, \boldsymbol{x_{u^{2}}}>,
\end{equation*}
\begin{equation*}
b_{11}= <\boldsymbol{x_{u^{1}u^{1}}}, \boldsymbol{n}>,\ \ \ b_{12}= <\boldsymbol{x_{u^{1}u^{2}}}, \boldsymbol{n}>,\ \ \ b_{22}= <\boldsymbol{x_{u^{2}u^{2}}}, \boldsymbol{n}>.
\end{equation*}
and $<,>$ is the Euclidean inner product.
For two sufficiently differentiable functions $p(u^{1}, u^{2})$ and $q(u^{1}, u^{2})$ on $M^{2}$ the first differential parameter of Beltrami with respect to the second fundamental form $II$ is defined by \cite{H1}
\begin{equation}  \label{6}
\nabla^{II}(p,q)=b^{ij}p_{/i}q_{/j},
\end{equation}
where $p_{/i}:=\frac{\partial p}{\partial u^{i}}$ and $b^{ij}$ are the components of the inverse tensor of $b_{ij}$. The second Beltrami operator according to the fundamental form $II$ of $M^{2}$ is defined by

\begin{equation}  \label{7}
\Delta ^{II}p =-b^{ij}\nabla^{II}_{i} p_{j}= -\frac{1}{\sqrt{|b|}}\frac{\partial}{\partial u^{i}}(\sqrt{|b|}b^{ij}\frac{\partial}{\partial u^{j}}),
\end{equation}
where $p$ is sufficiently differentiable function, $\nabla^{II}_{i}$ is the covariant derivative in the $u^{i}$ direction and $b= \det (b_{ij})$ \cite{A7}.

%

In the present paper, we mainly focus on surfaces of finite $II$-type by studying surfaces of revolution in $E^{3}$ which are connected, complete and of which their position vector $\boldsymbol{x}$ satisfies the following relation

\begin{equation}
\Delta ^{II}\boldsymbol{x}= A \boldsymbol{x},  \label{24}
\end{equation}
Our main result is
\begin{theorem} Spheres and catenoids are the only surfaces of revolution in $E^{3}$ whose position vector $\boldsymbol{x}$ satisfies condition (\ref{24}).
\end{theorem}
\section{Proof of the main theorem}

Let $C$ be a smooth curve lies on the $xz$-plane parametrized by
\begin{equation*}
\boldsymbol{r}(u) = (p(u), 0, q(u)), \ \ \ u\in (a, b),
\end{equation*}
where $p, q$ are smooth functions and $p$ is a positive function. When $C$ is revolved about the $z$-axis, the resulting point set $S$ is called the surface of revolution generated by the curve $C$. In this case, the $z$-axis is called the axis of revolution of $S$ and $C$ is called the profile curve of $S$. On the other hand, a subgroup of the rotation group which fixes the vector $(0, 0, 1)$ is generated by
\begin{equation*}
\left(
\begin{array}{ccc}
\cos v & -\sin v & 0 \\
\sin v & \cos v & 0 \\
0 & 0& 1
\end{array}%
\right).
\end{equation*}

Then the position vector of $S$ is given by
\begin{equation}  \label{3.1}
\boldsymbol{x}(u,v)= \big(p(u)\cos v, p(u)\sin v, q(u) \big), \ \ \ u\in (a,b), \ \ \ v \in [0,2\pi).
\end{equation}

(For the parametric representation of surfaces of revolution, see \cite{A11, D4, K8}).

Here, we may assume that $C$ has the arc-length parametrization, i.e., it satisfies

\begin{equation}  \label{3.2}
(p\prime)^{2}+ (q\prime)^{2}=1
\end{equation}
where $\prime := \frac{d}{du}$.
On the other hand $p\prime q\prime \neq 0$, because if $p$ = const. or $q$ = const. then $S$ is a circular cylinder or part of a plane, respectively. Hence the Gaussian curvature of $S$ vanishes. A case which has been excluded.

Using the natural frame $\{{\boldsymbol{x}_{u}, \boldsymbol{x}_{v}}\}$ of $S$ defined by
\begin{equation*}
\boldsymbol{x_{u}}=\left( p\prime(u)\cos v, p\prime(u)\sin v, q\prime(u)\right),
\end{equation*}
and
\begin{equation*} \newline
\boldsymbol{x_{v}}=\left( -p(u)\sin v, p(u)\cos v, 0\right),
\end{equation*}

\noindent the components $g_{ij}$ of the first fundamental form in (local) coordinates are the following
\begin{equation*}
g_{11}= 1,\ \ \ g_{12}=  0,\ \ \ g_{22}=  p^{2}.
\end{equation*}

Denoting by $R_{1}, R_{_{2}}$ the principal radii of curvature of $S$ and $\kappa$ the curvature of the curve $C$, we have
\begin{equation*}  \label{3.3}
R_{1} = \kappa, \ \ \ \ R_{2}= \frac{q\prime}{p}.
\end{equation*}

The mean and the Gaussian curvature of $S$ are respectively
\begin{align*}  \label{3.4}
2H = R_{1}+R_{2} = \kappa+ \frac{q\prime}{p},\ \ \ \ K = R_{1}R_{2} =\frac{\kappa q\prime}{p}= -\frac{p\prime\prime}{p}.
\end{align*}


The components $b_{ij}$ of the second fundamental form in (local) coordinates are the following
\begin{equation*}
b_{11} = \kappa,\ \ \ b_{12} = 0,\ \ \ b_{22} = pq\prime.
\end{equation*}

The Beltrami operator $\Delta ^{II}$ in terms of local coordinates $(u, v)$ of $S$ can be expressed as follows

\begin{eqnarray}  \label{3.5}
\Delta ^{II} =-\frac{1}{\kappa}\frac{\partial^{2}}{\partial u^{2}}-\frac{1}{pq\prime} \frac{\partial^{2}}{\partial v^{2}}
+\frac{1}{2}\bigg(\frac{\kappa\prime}{\kappa^{2}}-\frac{p\prime q\prime+\kappa pp\prime}{\kappa pq\prime}\bigg)\frac{\partial}{\partial u}.
\end{eqnarray}

On account of (\ref{3.2}) we can put
\begin{equation}  \label{3.6}
p\prime= \cos\varphi, \ \ \ q\prime = \sin\varphi,
\end{equation}
where $\varphi=\varphi(u)$. Then $\kappa = \varphi\prime$ and relation (\ref{3.5}) becomes
\begin{eqnarray}  \label{3.7}
\Delta ^{II} =-\frac{1}{\varphi\prime}\frac{\partial^{2}}{\partial u^{2}}-\frac{1}{p\sin\varphi} \frac{\partial^{2}}{\partial v^{2}}
+\frac{1}{2}\bigg(\frac{\varphi\prime\prime}{(\varphi\prime)^{2}}-\frac{\cos\varphi\sin\varphi+p\varphi\prime\cos\varphi}{p\varphi\prime\sin\varphi}\bigg) \frac{\partial}{\partial u},
\end{eqnarray}
while the mean and the Gaussian curvature of $S$ become
\begin{equation}   \label{3.8}
2H= \varphi\prime+\frac{\sin\varphi}{p},
\end{equation}

\begin{equation}   \label{3.9}
K= \frac{\varphi\prime\sin\varphi}{p}.
\end{equation}

Let $(x_{1}, x_{2}, x_{3})$ be the coordinate functions of $\boldsymbol{x}$ of (\ref{3.1}). Then we have
\begin{equation}  \label{3.10}
\Delta^{II}\boldsymbol{x} =(\Delta^{II}x_{1},\Delta^{II}x_{2},\Delta^{II}x_{3}).
\end{equation}

From (\ref{3.7}) and (\ref{3.10}), one can find

\begin{equation}  \label{3.11}
\Delta^{II}x_{1} =\Delta^{II}(p\cos v) = \bigg(\sin\varphi +\frac{1}{\sin\varphi}+\frac{\varphi\prime\prime\cos\varphi}{2\varphi\prime^{2}} - \frac{H\cos^{2}\varphi}{\varphi\prime\sin\varphi} \bigg)\cos v,
\end{equation}
\begin{equation}  \label{3.12}
\Delta^{II}x_{2} =\Delta^{II}(p\sin v) = \bigg(\sin\varphi +\frac{1}{\sin\varphi}+\frac{\varphi\prime\prime\cos\varphi}{2\varphi\prime^{2}} - \frac{H\cos^{2}\varphi}{\varphi\prime\sin\varphi} \bigg)\sin v,
\end{equation}
\begin{equation}  \label{3.13}
\Delta^{II}x_{3} =\Delta^{III}(q) = -\frac{3}{2}\cos\varphi +\frac{\varphi\prime\prime\sin\varphi}{2\varphi\prime^{2}} -\frac{\sin\varphi\cos\varphi}{2p\varphi\prime}.
\end{equation}

We denote by $a_{ij},i,j=1,2,3,$ the entries of the matrix $A.$ By using (\ref{3.11}), (\ref{3.12}) and (\ref{3.13}) condition (\ref{24}) is found to be equivalent to the following system
\begin{equation}  \label{3.14}
\bigg(\sin\varphi +\frac{1}{\sin\varphi}+\frac{\varphi\prime\prime\cos\varphi}{2\varphi\prime^{2}} - \frac{H\cos^{2}\varphi}{\varphi\prime\sin\varphi} \bigg)\cos v = a_{11}p\cos v +a_{12}p\sin v+ a_{13}q,
\end{equation}
\begin{equation}  \label{3.15}
\bigg(\sin\varphi +\frac{1}{\sin\varphi}+\frac{\varphi\prime\prime\cos\varphi}{2\varphi\prime^{2}} - \frac{H\cos^{2}\varphi}{\varphi\prime\sin\varphi} \bigg)\sin v = a_{21}p\cos v+a_{22}p\sin v+a_{23}q,
\end{equation}
\begin{equation}  \label{3.16}
-\frac{3}{2}\cos\varphi +\frac{\varphi\prime\prime\sin\varphi}{2\varphi\prime^{2}} -\frac{\sin\varphi\cos\varphi}{2p\varphi\prime} = a_{31}p\cos v+a_{32}p\sin v+a_{33}q.
\end{equation}

From (\ref{3.16}) it can be easily verified that $a_{31}= a_{32}= 0$. On differentiating (\ref{3.14}) and (\ref{3.15}) twice with respect to
$v$ we have
\begin{equation}  \label{3.17}
\bigg(\sin\varphi +\frac{1}{\sin\varphi}+\frac{\varphi\prime\prime\cos\varphi}{2\varphi\prime^{2}} - \frac{H\cos^{2}\varphi}{\varphi\prime\sin\varphi} \bigg)\cos v = a_{11}p\cos v +a_{12}p\sin v,
\end{equation}
\begin{equation}  \label{3.18}
\bigg(\sin\varphi +\frac{1}{\sin\varphi}+\frac{\varphi\prime\prime\cos\varphi}{2\varphi\prime^{2}} - \frac{H\cos^{2}\varphi}{\varphi\prime\sin\varphi} \bigg)\sin v = a_{21}p\cos v+a_{22}p\sin v.
\end{equation}

Thus $a_{13}q=a_{23}q=0,$ so that $a_{13}$ and $a_{23}$ vanish. Equations (\ref{3.14}), (\ref{3.15}) and (\ref{3.16}) are equivalent to the following
\begin{equation}  \label{3.19}
\bigg(\sin\varphi +\frac{1}{\sin\varphi}+\frac{\varphi\prime\prime\cos\varphi}{2\varphi\prime^{2}} - \frac{H\cos^{2}\varphi}{\varphi\prime\sin\varphi} \bigg)\cos v = a_{11}p\cos v +a_{12}p\sin v,
\end{equation}
\begin{equation}  \label{3.20}
\bigg(\sin\varphi +\frac{1}{\sin\varphi}+\frac{\varphi\prime\prime\cos\varphi}{2\varphi\prime^{2}} - \frac{H\cos^{2}\varphi}{\varphi\prime\sin\varphi} \bigg)\sin v = a_{21}p\cos v+a_{22}p\sin v,
\end{equation}
\begin{equation}  \label{3.21}
-\frac{3}{2}\cos\varphi +\frac{\varphi\prime\prime\sin\varphi}{2\varphi\prime^{2}} -\frac{\sin\varphi\cos\varphi}{2p\varphi\prime} = a_{33}q.
\end{equation}

But $\sin v$ and $\cos v$ are linearly independent functions of $v$, so we finally obtain $a_{12}=a_{21}=0,a_{11}=a_{22}$. Putting $a_{11}=a_{22}=\lambda$ and $a_{33}=\mu$, we see that the system of equations (\ref{3.19}), (\ref{3.20}) and (\ref{3.21}) reduces now to the following two equations

\begin{equation}  \label{3.22}
\sin\varphi+\frac{1}{\sin\varphi}+\frac{\varphi\prime\prime\cos\varphi}{2\varphi\prime^{2}}-\frac{H\cos^{2}\varphi}{\varphi\prime\sin\varphi}= \lambda p,
\end{equation}
\begin{equation}  \label{3.23}
-\frac{3}{2}\cos\varphi +\frac{\varphi\prime\prime\sin\varphi}{2\varphi\prime^{2}} -\frac{\sin\varphi\cos\varphi}{2p\varphi\prime} = \mu q.
\end{equation}

Hence the matrix $A$ for which relation (\ref{24}) is satisfied becomes
\begin{equation*}
A =\left[
\begin{array}{ccc}
\lambda & 0 & 0 \\
0 & \lambda & 0 \\
0 & 0& \mu
\end{array}%
\right].
\end{equation*}

We distinguish the following cases:

\textit{Case I.}  $\lambda=\mu =0$.
Equations (\ref{3.22}) and (\ref{3.23}) become
\begin{equation}  \label{3.24}
\sin\varphi +\frac{1}{\sin\varphi}+\frac{\varphi\prime\prime\cos\varphi}{2\varphi\prime^{2}} - \frac{H\cos^{2}\varphi}{\varphi\prime\sin\varphi}= 0,
\end{equation}
\begin{equation}  \label{3.25}
-\frac{3}{2}\cos\varphi +\frac{\varphi\prime\prime\sin\varphi}{2\varphi\prime^{2}} -\frac{\sin\varphi\cos\varphi}{2p\varphi\prime} = 0.
\end{equation}

Multiplying (\ref{3.24}) by $\sin\varphi$, and (\ref{3.25}) by $-\cos\varphi$ and adding the resulting of these two equations whence it follows $\cos^{2}\varphi +\sin^{2}\varphi + 1= 0$, a contradiction.

\textit{Case II.}  $\lambda=\mu \neq0$.
Equations (\ref{3.22}) and (\ref{3.23}) become
\begin{equation}  \label{3.26}
\sin\varphi+\frac{1}{\sin\varphi}+\frac{\varphi\prime\prime\cos\varphi}{2\varphi\prime^{2}}-\frac{H\cos^{2}\varphi}{\varphi\prime\sin\varphi}=\lambda p,
\end{equation}
\begin{equation}  \label{3.27}
-\frac{3}{2}\cos\varphi +\frac{\varphi\prime\prime\sin\varphi}{2\varphi\prime^{2}} -\frac{\sin\varphi\cos\varphi}{2p\varphi\prime} = \lambda q.
\end{equation}

Similarly, multiplying (\ref{3.26}) by $\sin\varphi$, and (\ref{3.27}) by $-\cos\varphi$ and adding the resulting of these two equations whence it follows
\begin{equation}  \label{3.28}
\lambda p\sin\varphi -\lambda q\cos\varphi = 2.
\end{equation}

On differentiating the last equation with respect to $u$ we find
\begin{equation}  \label{3.29}
\lambda\varphi\prime(pp\prime + qq\prime) = 0.
\end{equation}

$\lambda\varphi\prime$ cannot be equal 0 otherwise, from (\ref{3.9}) the Gauss curvature vanishes. Hence, $pp\prime + qq\prime=0$, i.e. $(p^{2}+q^{2})\prime =0$. Therefore $(p^{2}+q^{2}) =const.$. Thus $C$ is part of a circle and $S$ is obviously part of a sphere.

\textit{Case III.}  $\lambda\neq0, \mu =0$.
Following the same procedure as in $Case I$ and $Case II$, we can obtain

\begin{equation}  \label{3.30}
2-\lambda p\sin\varphi=0.
\end{equation}

Differentiating (\ref{3.30}) with respect to $u$ we have
\begin{equation}  \label{3.31}
(\sin\varphi+p\varphi\prime)\cos\varphi = 0.
\end{equation}

Taking into account relation (\ref{3.8}), equation (\ref{3.31}) becomes
\begin{equation}  \label{3.32}
2Hp\cos\varphi = 0
\end{equation}
which implies the mean curvature $H$ vanishes identically. Therefore, the surface is minimal, that is, it is a catenoid. Furthermore, a
catenoid satisfies the condition (\ref{24}).

\textit{Case IV.}  $\lambda=0, \mu \neq0$.
In this case (\ref{3.22}) and (\ref{3.23}) are given respectively by
\begin{equation}  \label{3.33}
\sin\varphi +\frac{1}{\sin\varphi}+\frac{\varphi\prime\prime\cos\varphi}{2\varphi\prime^{2}}- \frac{H\cos^{2}\varphi}{\varphi\prime\sin\varphi}= 0,
\end{equation}
\begin{equation}  \label{3.34}
-\frac{3}{2}\cos\varphi +\frac{\varphi\prime\prime\sin\varphi}{2\varphi\prime^{2}} -\frac{\sin\varphi\cos\varphi}{2p\varphi\prime} = \mu q.
\end{equation}

Following the same procedure as in $Case I$ and $Case II$, we find

\begin{equation}  \label{3.35}
2+\mu q\cos\varphi=0.
\end{equation}

Differentiating this equation we have
\begin{equation}  \label{3.36}
q\varphi\prime- \cos\varphi=0,
\end{equation}
from which
\begin{equation}  \label{3.37}
\varphi\prime=\frac{\cos\varphi}{q}.
\end{equation}

Another differentiation of (\ref{3.36}), gives
\begin{equation}  \label{3.38}
2\varphi\prime\sin\varphi+q\varphi\prime\prime=0.
\end{equation}

From (\ref{3.37}) and (\ref{3.38}), we have
\begin{equation}  \label{3.39}
\varphi\prime\prime=-\frac{2\sin\varphi\cos\varphi}{q^{2}}.
\end{equation}

Equation (\ref{3.33}) can be written
\begin{equation}  \label{3.40}
1+\sin^{2}\varphi +\frac{\varphi\prime\prime\cos\varphi\sin\varphi}{2\varphi\prime^{2}}-\frac{1}{2}\cos^{2}\varphi -\frac{\cos^{2}\varphi\sin\varphi} {2p\varphi\prime}= 0.
\end{equation}

Consequently, from (\ref{3.35}), (\ref{3.37}) and (\ref{3.39}), equation (\ref{3.40}) becomes
\begin{equation}  \label{3.41}
2-\cos^{2}\varphi+\frac{2\sin\varphi}{\mu p} = 0,
\end{equation}
from which
\begin{equation}  \label{3.42}
\mu p=\frac{2\sin\varphi}{\cos^{2}\varphi -2}.
\end{equation}

Differentiating (\ref{3.42}), we get
\begin{equation}  \label{3.43}
\mu =\frac{2\varphi\prime}{\cos^{2}\varphi -2}+\frac{4\varphi\prime\sin^{2}\varphi}{(\cos^{2}\varphi -2)^{2}}.
\end{equation}

On using (\ref{3.37}) and (\ref{3.35}) after some computation, we can obtain that $\sin\varphi=0$, that is, $q = const.$, which implies that the Gauss curvature vanishes. A case which was excluded. Thus, there are no surfaces of revolution satisfying this case.

\textit{Case V.}  $\lambda\neq0, \mu \neq0$.
If we multiply (\ref{3.26}) by $\sin\varphi$, and (\ref{3.27}) by $-\cos\varphi$ and adding the resulting equations, we easily get
\begin{equation}  \label{3.44}
\lambda p\sin\varphi -\mu q\cos\varphi = 2.
\end{equation}

We put
\begin{equation}  \label{3.45}
\Omega:=\lambda p\sin\varphi +\mu q\cos\varphi .
\end{equation}

By using (\ref{3.44}), the derivative of $\Omega$ is the following
\begin{equation}  \label{3.46}
\Omega\prime=\lambda \cos^{2}\varphi +\mu \sin^{2}\varphi -2\varphi\prime.
\end{equation}

Differentiating the equation (\ref{3.44}) and using (\ref{3.45}) we find
\begin{equation}  \label{3.47}
\Omega\varphi\prime=(\mu -\lambda)\cos\varphi \sin\varphi.
\end{equation}

It is easily verified that $\Omega \neq0$, hence (\ref{3.47}) can be written
\begin{equation}  \label{3.48}
\varphi\prime=\frac{(\mu -\lambda)\cos\varphi \sin\varphi}{\Omega}.
\end{equation}

Differentiating the last equation and using (\ref{3.46}) and (\ref{3.47}) we obtain
\begin{equation}  \label{3.49}
\varphi\prime\prime=\frac{\big((\lambda -2\mu)\sin^{2}\varphi+(\mu -2\lambda)\cos^{2}\varphi\big)\varphi\prime+2\varphi\prime^{2}}{\Omega}.
\end{equation}

In view of (\ref{3.48}) and (\ref{3.49}) relation (\ref{3.23}) takes the following form
\begin{equation}  \label{3.50}
\frac{\Omega\cos\varphi}{p}+\frac{2(\lambda -\mu)\cos\varphi \sin\varphi}{\Omega}-2\mu(\lambda -\mu)q\cos\varphi -(\lambda -2\mu)=0,
\end{equation}

Multiplying the last equation by $\Omega p\cos\varphi $, we have
\begin{eqnarray}  \label{3.52}
&&2(\lambda -\mu)p\cos^{2}\varphi \sin\varphi-2\mu(\lambda-\mu)pq\Omega\cos^{2}\varphi  \notag \\
&&-(\lambda -2\mu)\Omega p\cos\varphi+ \Omega^{2}\cos^{2}\varphi =0.
\end{eqnarray}

From (\ref{3.44}) and (\ref{3.45}), it can be is easily verified that
\begin{equation}  \label{3.53}
\Omega\cos\varphi =\lambda p-2\sin\varphi.
\end{equation}

Therefore, on using (\ref{3.44}) and (\ref{3.53}), relation (\ref{3.52}) becomes
\begin{equation}  \label{3.54}
a_{1}p^{3}+ a_{2}p^{2}+a_{3}p+a_{4}=0,
\end{equation}
where, here we put
\begin{equation*}
a_{1}= \lambda^{2}(\mu - \lambda)\sin\varphi, \ \ \ \ \ \ a_{2}=\lambda[(2\lambda-\mu)-2(\mu - \lambda)\sin^{2}\varphi],
\end{equation*}
\begin{equation*}
a_{3}= [(\mu - \lambda)\sin^{2}\varphi-(\mu - 4\lambda)]\sin\varphi, \ \ \ \ \ \ a_{4}=2\sin^{2}\varphi.
\end{equation*}

Taking the derivative of (\ref{3.54}) and then by using (\ref{3.45}), (\ref{3.48}) and (\ref{3.53}), we obtain
\begin{equation}  \label{3.55}
b_{1}p^{3}+ b_{2}p^{2}+b_{3}p+b_{4}=0,
\end{equation}
where
\begin{equation*}
b_{1}= \lambda^{2}(\mu - \lambda)\sin\varphi[(2\lambda+\mu)-(\mu - \lambda)\sin^{2}\varphi],
\end{equation*}
\begin{equation*}
b_{2}= 2\lambda[\lambda(2\lambda-\mu)-(\mu - \lambda)(3\lambda+2\mu)\sin^{2}\varphi+2(\mu - \lambda)^{2}\sin^{4}\varphi],
\end{equation*}
\begin{equation*}
b_{3}= [(8\lambda\mu - 8\lambda^{2}-\mu^{2})+2(\mu - \lambda)(2\mu +\lambda)\sin^{2}\varphi-3(\mu - \lambda)^{2}\sin^{4}\varphi]\sin\varphi,
\end{equation*}
\begin{equation*}
b_{4}= 6(\mu - 2\lambda)\sin^{2}\varphi-6(\mu - \lambda)\sin^{4}\varphi.
\end{equation*}

Combining (\ref{3.54}) and (\ref{3.55}) we conclude that
\begin{equation}  \label{3.56}
c_{1}p^{2}+c_{2}p+c_{3}=0,
\end{equation}
where
\begin{equation}  \label{3.57}
c_{1}= a_{1}b_{2}-a_{2}b_{1}=\lambda[2(\mu - \lambda)^{2}\sin^{4}\varphi-3\mu(\mu - \lambda)\sin^{2}\varphi-\mu(2\lambda-\mu)],
\end{equation}
\begin{eqnarray}  \label{3.58}
c_{2}&=&a_{1}b_{3}-a_{3}b_{1}=2[-(\mu - \lambda)^{2}\sin^{4}\varphi   \notag \\
&&+(\mu+2\lambda)(\mu - \lambda)\sin^{2}\varphi+\lambda(3\mu -8\lambda)]\sin\varphi,
\end{eqnarray}
\begin{equation}  \label{3.59}
c_{3}=a_{1}b_{4}-a_{4}b_{1}=[-4(\mu - \lambda)\sin^{4}\varphi+4(\mu - 4\lambda)\sin^{2}\varphi].
\end{equation}

Taking the derivative of (\ref{3.56}) and then by using (\ref{3.45}), (\ref{3.48}) and (\ref{3.53}), we obtain
\begin{equation}  \label{3.60}
d_{1}p^{2}+ d_{2}p+d_{3}=0,
\end{equation}
where
\begin{equation*}
d_{1}= -4\lambda(\mu - \lambda)^{3}\sin^{6}\varphi +\sum^{2}_{i=0}D_{1i}(\lambda,\mu)\sin^{2i}\varphi,
\end{equation*}
\begin{equation*}
d_{2}= 5(\mu - \lambda)^{3}\sin^{7}\varphi +\sum^{2}_{i=0}D_{2i}(\lambda,\mu)\sin^{2i+1}\varphi,
\end{equation*}
\begin{equation*}
d_{3}= 10(\mu - \lambda)^{2}\sin^{6}\varphi +\sum^{2}_{i=0}D_{3i}(\lambda,\mu)\sin^{2i}\varphi,
\end{equation*}
and $D_{ji}(\lambda,\mu),\ (j = 1, 2,3)$ are polynomials in $\lambda$ and $\mu$. Combining (\ref{3.56}) and (\ref{3.60}) we find that
\begin{equation}  \label{3.61}
e_{1}p+ e_{2}=0,
\end{equation}
where
\begin{equation}  \label{3.62}
e_{1}= c_{1}d_{2}-c_{2}d_{1}=2(\mu - \lambda)^{5}\sin^{10}\varphi+\sum^{4}_{i=0}E_{1i}(\lambda,\mu)\sin^{2i}\varphi,
\end{equation}
\begin{equation}  \label{3.63}
e_{2}=c_{1}d_{3}-c_{3}d_{1}=20(\mu - \lambda)^{4}\sin^{9}\varphi +\sum^{3}_{i=0}E_{2i}(\lambda,\mu)\sin^{2i+1}\varphi,
\end{equation}
and $E_{ji}(\lambda,\mu),\ (j = 1, 2)$ are some polynomials in $\lambda$ and $\mu$.
Following the same procedure by taking the derivative of (\ref{3.61}) and taking into account (\ref{3.45}), (\ref{3.48}) and (\ref{3.53}), we find
\begin{equation}  \label{3.64}
h_{1}p+ h_{2}=0,
\end{equation}
where
\begin{equation}  \label{3.65}
h_{1}=-20(\mu - \lambda)^{6}\sin^{12}\varphi+\sum^{5}_{i=0}H_{1i}(\lambda,\mu)\sin^{2i}\varphi,
\end{equation}
\begin{equation}  \label{3.66}
h_{2}=-184(\mu - \lambda)^{5}\sin^{11}\varphi +\sum^{4}_{i=0}H_{2i}(\lambda,\mu)\sin^{2i+1}\varphi,
\end{equation}
and $H_{ji}(\lambda,\mu),\ (j = 1, 2)$ are polynomials in $\lambda$ and $\mu$. Combining (\ref{3.61}) and (\ref{3.64}) we finally find
\begin{equation}  \label{3.67}
32(\mu - \lambda)^{10}\sin^{20}\varphi+\sum^{9}_{i=0}P_{i}(\lambda,\mu)\sin^{2i}\varphi =0.
\end{equation}
where $P_{i}(\lambda,\mu),\ (i = 0, 1, ...,9)$  are the known polynomials in $\lambda$ and $\mu$. Since this polynomial is equal to zero for every $\varphi$, all its coefficients must be zero. Therefore, we conclude that $\mu - \lambda = 0$, which is a contradiction. Consequently, there are no surfaces of revolution in this case. This completes our proof.






\end{document}